\documentclass{primus}

\usepackage{natbib}
\usepackage[hyphens,spaces]{url}
\usepackage[font=footnotesize,labelfont=bf]{caption}
\usepackage{amsmath}
\usepackage{url}
\usepackage{amssymb,amsthm}
\usepackage{graphicx}
\usepackage{enumitem}
\usepackage{caption}
\usepackage[margin=3cm]{geometry}
\usepackage{xcolor}

\def\Z{\mathbb{Z}}
\newcommand{\N}{\mathbb{N}}

\newtheorem{theorem}{Theorem}
\theoremstyle{definition}
\newtheorem{example}[theorem]{Example}

\begin{document}

\title{Evaluating the use of e-assessment in a first-year pure mathematics module}
\author{Stefanie Zegowitz}
\maketitle

\begin{abstract}
\noindent 
This article presents the findings of a case study which introduced online quizzes as a form of assessment in pure mathematics. Rather than being designed as an assessment \textit{of} learning, these quizzes were designed to be an assessment \textit{for} learning; they aimed to academically support students in their transition from A-Level mathematics to university-level pure mathematics by providing an extrinsic motivation to engage them with their learning material early on and to emphasize the small details within proofs, such as defining notation, which are not necessarily emphasized by written homework assignments. The results obtained during the two-year study using online quizzes show e-assessment to be a powerful complementary tool to traditional written homework assignments. 
\end{abstract}

\listkeywords{e-assessment, mathematical proofs, undergraduate mathematics education}

\section{Introduction}

Mathematics is a discipline with unique features when it comes to assessment and teaching \cite[]{lms}. It employs few summative assessment methods at university-level with the main method being timed closed book examinations \cite[]{paola}. Researchers have often called for a larger variety of its summative assessment methods \cite[]{steen}, and \cite{birenbaum} suggest that traditional timed exams do not reflect the future needs of today's learners as they seem to be designed as assessments \textit{of} learning rather than \textit{for} learning.

E-assessments offer several advantages, especially when teaching large cohorts: they offer the possibility of automatic marking and feedback. However, the majority of e-learning systems are poorly adapted for the use in mathematics, a language in its own right \cite[]{gruttmann}. According to \cite{smith}, current e-learning systems do not adequately support the necessary notation and diagrams, ``the very building blocks of math communication" (p.~324).

In this paper, I will present the findings from a case study investigating the use of online quizzes in pure mathematics. These bi-weekly open book online quizzes formed part of the summative assessment of a first-year undergraduate mathematics module in the UK, counting towards 10\% of the final grade. The quizzes provided automatic and individual feedback, tailored to selected answer choices, making the quizzes summative as well as formative assessment tools designed \textit{for} learning.

The case study is based on a teaching intervention I designed and implemented as a lecturer in one of my modules, Introduction to Proofs, which introduced mathematics students to the new language of pure mathematics. The case study has two motivations which are based on a couple of observations I made when teaching the module. Firstly, I observed that students appeared to be initially overwhelmed by the transition from A-Level mathematics to university-level pure mathematics: they appeared worried and stressed, repeatedly mentioning how different pure mathematics is compared to A-level mathematics. According to \cite{anthony}, the study for tests and exams is a perceived key factor, by lecturers and students, for the successful transition to university-level mathematics. Therefore, one aim of the case study was to provide students with an extrinsic motivation to engage with learning material early on by changing the assessment method of the module -- which was previously solely assessed via a summative end-of-term exam -- to additionally introduce periodic online quizzes. 

Secondly, I observed that written formative homework did not necessarily emphasize the importance of the small details within proofs. A classic example of this is the proof of $ \sqrt{2} $ being irrational. Here, we define $ \sqrt{2}=\frac{a}{b} $ where $ a\in\Z$ and $ b\in\Z\backslash\{0\} $ are relatively prime, and we continue to argue that both $ a $ and $ b $ are even which contradicts our initial assumption that $ a $ and $ b $ are relatively prime. Had we, for example, not defined $ a$ and $b$ to be relatively prime, then there would not have been a contradiction and the whole argument would have collapsed. Henceforth, another motivation of the case study was to emphasize the small details within proofs, such as defining notation.

\section{Background}

E-assessments have two key advantages: they substantially reduce marking time, and secondly, they allow for immediate feedback. Feedback supports optimal learning and therefore plays a critical role in e-assessment \cite[]{gaytan}. The integration of formative assessment, such as feedback, presents an opportunity to provide continuous student support, thereby nurturing meaningful engagement and deep learning \cite[]{gikandi}. \cite{wolsey} further shows that its use in an online setting supports students to identify their strengths and weaknesses, revise their work, and continuously refine their understanding.

Several studies have shown that e-assessments have a positive effect on summative exam scores and overall academic success \cite[]{angus, salas, dobson, kibble}. This effect seems to be emphasized when e-assessments are scheduled periodically and when they only count towards a small percentage of the overall mark. Periodic scheduling is especially beneficial since it continuously and consistently engages students with their learning material \cite[]{savander}. In addition, it has the benefit to encourage distributive learning as opposed to cramming \cite[]{desouza}. 

Focusing on pure mathematics, we find that its foundation is built on mathematical proofs. Conclusively, the main method to assess a mathematics student's performance at university-level is to assess their ability to construct proofs \cite[]{weber}. The predicament with assessing proofs is that mathematicians need to not only evaluate how well a particular proof is written but also how well it is understood \cite[]{miller}, and there is a definite lack of assessment tools on proof comprehension \cite[]{ramos}. How we assess a student's proof comprehension remains an open question \cite[]{ramos}, but what we do know is that feedback seems to be one of the main tools which supports students in learning how to construct valid proofs \cite[]{moore}. Proof validation, which is the process of reading and reflecting upon proofs to determine their correctness \cite[]{selden}, is another tool which has been shown to have a positive impacts on students' ability to construct their own proofs \cite[]{powers}. It requires students to engage thoroughly with their learning material by, for example, asking and answering questions, constructing sub-proofs and finding and interpreting definitions and theorems \cite[]{selden}. 

While the main method of assessing a student's ability to construct proofs is via timed closed books examinations, an unusual way to assess a student's ability to validate proofs is via multiple-choice quizzes. This type of assessment emphasises the importance of rigour in proofs and allows for rapid and detailed feedback. \cite{hetherington} shows that the provision of feedback in paper based multiple-choice tests makes students not only be aware of what answer choice is correct but also be aware of common misconceptions in proofs. One may easily conclude that this type of assessment would be rather suited for e-learning systems. However, the general literature seems to agree that the majority of e-learning systems is rather poorly adapted to mathematics and its notation \cite[]{gruttmann, hetherington, smith}. Therefore, in mathematics, multiple-choice questions cannot be easily transferred from paper to an e-learning system. In addition, we have the predicament that available e-learning systems are not able to automatically assess a student's proof comprehension \cite[]{hetherington, sangwin}. So despite their benefits, e-learning systems are rarely used for the assessment of proofs. 

As a matter of fact, only one system exists which is dedicated to the assessment of proofs in pure mathematics, the system EASy. Specifically, EASy is designed to assess proving skills by supervising the correct application of proof strategies and giving feedback when rules are incorrectly applied \cite[]{gruttmann}. While it is thought to support the learning of elementary proving skills, EASy also seems to have some drawbacks: It is slightly complex to use and prolongs solving time in comparison to solving a proof manually \cite[]{gruttmann}. In subsequent developments, EASy seems to have moved away from the assessment of proofs towards the assessment of computer programming \cite[]{sangwin}. Other e-learning systems in pure mathematics do exist, for example Theorema \cite[]{buchberger} which ``provides the mathematical language and the possibility of fully automated generation of mathematical proofs'' \cite[][p.~1]{mayrhofer}, but these systems tend to not focus on assessment.

Another system worth mentioning is STACK, a \textit{System for Teaching and Assessment using a Computer algebra Kernel}, which can be applied to a wide range of mathematics. STACK allows students to submit answers in the form of algebraic expressions. Once an expression is entered, the system then establishes its properties and provides feedback, a numerical mark and an “internal ‘note’ for later analysis” \cite[][p.~103]{sangwin}. Using STACK to turn mathematical proofs into interactive online quizzes, \cite{siri} show that students “understood the proofs much more clearly than when the proofs were recited during lectures” and found them to be “as useful as doing traditional problem sheets” (p.~9).

While e-learning system are still rarely used for the assessment of proofs, \cite{greenhow} nevertheless recommends the use of an e-learning system ``as part of a blended formative and summative assessment of any mathematics module at level 1 of a university degree" (p.~134). In order to successfully introduce such a system, he recommends the encoding of common students' misconceptions, the use of random parameters and the giving of detailed feedback which not only presents a model solution but also provides individual feedback tailored to selected answer choices. 

This article investigates the use of online quizzes as a form of assessment \textit{for} learning in pure mathematics. It evaluates the following research questions:
\begin{itemize}
\item How does periodic e-assessment affect student engagement with learning material?
\item Does periodic e-assessment help to emphasize the small details within proofs, such as defining notation?
\item Consequently, does periodic e-assessment help students with the subsequent writing of their own proofs? 
\end{itemize}
The following assumptions are made, based on constructivist theory:
\begin{itemize}
\item[$ \circ $]
Learners construct new knowledge upon the foundation of their previous knowledge \cite[]{hoover}
\item[$ \circ $]
Learning occurs through learners' active involvement in the construction of meaning \cite[]{piaget}
\end{itemize}
Note that here we interpret `active involvement in the construction of meaning' to be `engagement with learning material'. More precisely, by `engagement' we mean `interaction'. These interactions include, for example, critical reading.

\section{Implementation}
\label{implementation}

The teaching intervention was implemented in a first-year UK undergraduate mathematics module, Introduction to Proofs, which I taught for three years. This module introduced a large cohort of approximately 300 to 350 first-year undergraduate mathematics students to the language of pure mathematics, and it had an emphasis on definitions, proofs and their techniques. It covered fundamental concepts such as logical propositions, basic set theory and cardinality, functions and relations, and topics from elementary number theory. 

Initially, the module was assessed solely via a summative end-of-term exam and formative weekly written homework assignments, which were encouraged but not compulsory. In my second year of teaching, I introduced bi-weekly online quizzes, which formed part of the module's summative assessment and which replaced every other written homework assignment.

According to \cite{heijne}, open book assessments enable students to demonstrate a better understanding of learning material since they are able to answer questions ``by reading and thinking rather than reading and memorizing" (p.~16). Other studies have shown that open book assessments reduce the need for cramming and thereby cause less stress and anxiety before and during an assessment \cite[]{feldhusen, michaels, theo}. Since one of the aims of the online quizzes was for students to meaningfully engage with their learning material as opposed to students cramming and merely memorising their learning material, these studies supported my choice for the online quizzes to be of open book form. 

Being the university's choice of e-learning platform for students, I used Blackboard to implement the quizzes. While the assessment system STACK is dedicated to the e-learning platform Moodle, working with Blackboard nevertheless had several advantages -- it allowed the use of \LaTeX, automatic marking, automatic and individual feedback tailored to students' selected answer choices, and it supported a variety of question types. In addition, it offered the feature to randomise in which order questions and answer choices appeared. This had the advantage that each students would see a different version of the same quiz. 

For each quiz, I set a total of 7 questions: 4 easier questions counting 2 points each and 3 harder questions counting 4 points each, adding up to a total of 20 points per quiz. The easier questions aimed to familiarise students with definitions and to test knowledge, therefore being a local type of assessment according to the assessment model for proof comprehension proposed by \cite{ramos}. The harder questions aimed to emphasize the small details and structures of mathematical proofs by asking students to either validate proofs or to arrange a proof in order. According to the assessment model by \cite{ramos}, these questions were a combination of local as well as holistical assessment. To construct the easier questions, I used my experience on the module to identify commonly misunderstood definitions and results. To construct the harder questions, I identified students' common misconceptions when writing mathematical proofs. Some of these include the following:

\begin{itemize}
\item[$ \circ $]
not defining notation, for example, not defining that the element $a$ belongs to the set $A$ or that $b$ is an integer
\item[$ \circ $]
working backwards in a proof:  assuming what needs to be proved, showing that this assumption leads to a true statement and then erroneously concluding that this means that the initial assumption must have been correct
\item[$ \circ $]
incorrect assumptions in a proof by contradiction so that, in fact, there is no contradiction
\item[$ \circ $]
only proving one half of an `if and only if' statement
\item[$ \circ $]
incorrectly forming the contrapositive
\item[$ \circ $]
not covering all cases in a proof
\item[$ \circ $]
using the wrong base case in a proof by induction
\item[$ \circ $]
using regular induction instead of strong induction
\end{itemize}

To be able to address most items in the above list, I used different types of questions: multiple-choice, fill-in-the-blank, matching and multiple-answer questions. Multiple-choice questions are great for emphasising the small details within a proof, such as defining notation or using correct initial assumption. Since an answer choice is either correct or incorrect, this question emphasises how an entire proof may collapse due to one small detail. An example of a multiple-choice question is below. For this example, the correct answer choice is c). Please refer to Example \ref{feedback} for the corresponding feedback.

\begin{example}
\label{relations}
Let $ a,b\in\Z. $ Define a relation on $\Z$ by $x\sim y$ if and only if $8x\equiv 5y $ (mod $3).$ Which one of the following proofs is correct? 

\begin{enumerate}[label=\alph*)]
\item
We will show that $\sim$ is symmetric. 

Let $x,y\in \Z$ such that $ x\sim y. $ Then $8x\equiv 5y $ (mod $3).$ Since $ 8 \equiv -1 $ (mod $3),$ then  $ 8x\equiv -x $ (mod $3$). Since $ 5 \equiv -1 $ (mod $3),$ then  $ 5y\equiv -y $ (mod $3$). Hence, we have that
\[
-x+y\equiv 8x-5y\equiv 0 \text{ (mod } 3),
\]
so $3\mid -x+y.$ Hence, there exists $n\in\mathbb{N}$ such that $-x+y=3n.$ Then
\[
8y-5x=8(3n+x)-5x=24n-3x=3(8n-x).
\]
Since $8n-x\in\Z,$ then $ 3\mid 8y-5x $ which means that $ 8y\equiv 5x $ (mod $3$). Hence, $ y\sim x, $ so $ \sim $ is symmetric.
\item
We will show that $\sim$ is reflexive. 

We have that $ 8x-5x=3x. $ Hence, $ 3\mid 8x-5x, $ so $ 8x\equiv 5x $ (mod $3$). This means that $x\sim x,$ so $\sim$ is reflexive.
\item 
We will show that $\sim$ is transitive. 

Let $ x,y,z\in\Z $ such that $ x\sim y $ and $ y\sim z. $ Then $ 8x\equiv 5y $ (mod $3$) and $ 8y\equiv 5z $ (mod $3$). This means that $ 3\mid 8x-5y $ and $ 3\mid 8y-5z,  $ so there exist $ m,n\in\Z $ such that $ 8x-5y=3m $ and $ 8y-5z=3n. $ Then
\[
8x-5z= (3m+5y) +(3n-8y)=3m+3n-3y=3(m+n-y).
\] 
Since $ m+n-y\in\Z, $ then $ 3\mid 8x-5z. $ This means that $ 8x\equiv 5z $ (mod $3$), so $ x\sim z. $ Therefore, $ \sim $ is transitive.
\item 
We will show that $ \sim $ is symmetric.

Let $ x\in\Z. $ Since $ 8\equiv 5 $ (mod $3$), then $ 8x\equiv 5x $ (mod $3$). Hence, we have that $ x\sim x, $ so $ \sim $ is symmetric.
\end{enumerate}
\end{example}

While multiple-choice questions are also great for asking students to correctly identify a definition or a result, fill-in-the-blank questions are great for asking students to apply such a definition or a result. This question type diminishes the effect of guessing that other types may allow. An example of a fill-in-the-blank question is below.

\begin{example}
\label{fillintheblank}
Please fill in the following blank: Let $ x\in\Z $ with $ -10\le x\le 50. $ If $ x\equiv 3 $ (mod$5$) and $ x\equiv 2 $ (mod $7$), then $ x= \rule{1cm}{0.15mm}\,. $
\end{example}

Matching questions support the setting of several questions, for example, about a particular set or a function. On the other hand, they allow for several statements to be matched with certain properties such as the property of a function being bijective or not. An example of a matching question is below.

\begin{example}
\label{matching1}
Let $ \{A_i: i\in\N\} $ and $ \{B_i:i\in\N\} $ be a collection of subsets of $ \N $ where we define $ A_i=\{2z+2: z\in\Z, z\ge \frac{i-2}{2}\} $ and $ B_i=\{2z+3:z\in\Z, z\ge \frac{i-3}{2}\}. $ Please match the following equations:
\begin{flushleft}
\begin{tabular}{lllll}
1. & $ \bigcup\limits_{i\in\N} A_i= $ & & a) & $ \N $ \\
2. & $ \bigcap\limits_{i\in\N} B_i= $ & & b) & $ 2\N $ \\
3. & $ \bigcup\limits_{i\in\N} (A_i\cup B_i)= $ & & c) & $ 2\N-1 $ \\
  &  & & d) & $ \emptyset $
\end{tabular}
\end{flushleft}
\end{example}

Matching questions are also great for asking students to arrange a proof in order. Arranging a proof in order emphasises structural issues such as working backwards or using incorrect initial assumptions in a proof by contradiction. For this type of question, students were given one and only one statement concerning what the proof will show. An example of such a matching question is below. More examples of different question types can be found in the appendix.

\begin{example}
\label{matching2}
Define $ f:\N\to\Z $ by $ f(x)=x^2-4. $ Please arrange the following proof by \textbf{contradiction} in the correct order.
\begin{flushleft}
\begin{tabular}{p{1cm} p{1cm} p{0.5cm} p{11.5cm}}
Step 1 & & a) & Therefore, $ f(x)\neq f(y). $ \\
Step 2 & & b) & Take $ x,y\in\N $ such that $ x=y. $ \\
Step 3 & & c) & We will show that $ f $ is injective. \\
Step 4 & & d) & Then $ x^2-4=f(x)=f(y)=x^2 $ which holds if and only if \\
& & & $ x^2-4=y^2-4 $ which holds if and only if $ x^2=y^2. $ \\
Step 5 & & e) & Then $ x=\pm y. $ \\
Step 6 & & f) & It follows that $ f $ is injective. \\
Step 7 & & g) & Take $ x,y\in\N $ such that $ x\neq y. $ \\
Step 8 & & h) & To the contrary, suppose that $ f(x)=f(y). $ \\
& & i) & This is a contradiction since $ x\neq y. $ \\
& & j) & Then $ x^2-4=f(x)\neq f(y)=x^2 $ which holds if and only if \\
& & & $ x^2-4\neq y^2-4 $ which holds if and only if $ x^2\neq y^2. $ \\
& & k) & To the contrary, suppose that $ f(x)\neq f(y). $ \\
& & l) & Since $ x,y\in\N, $ then $ x=y. $ \\
& & m) & Therefore, $ f(x)\neq f(y). $ \\
& & n) & Since $ x,y\in\N, $ then $ x\neq y. $ 
\end{tabular}
\end{flushleft}
\end{example}

Initially, students were given 48 hours to complete each quiz. They were allowed to edit and save their answers as many times as they wished, but they were only allowed to submit a quiz once. In my second year of using online quizzes, the deadline was extended to one week in consideration of student feedback and students' concurrent coursework deadlines.

Since the online quizzes were bi-weekly, they would typically address two weeks worth of material covered in lectures.  According to \cite{zhu}, their work provides substantial evidence that students who learn via sequences of mathematical examples are ``at least as successful as, and sometimes more successful than, students learning by conventional methods" (p.~160). They further claim that learning via examples is ``learning with understanding" which enables students to recognise what rules to apply to a particular example and to apply that rule (p.~160). Therefore, during weekly problem classes, I included examples similar to questions which would appear in the upcoming quiz. 

In addition to learning from examples, students also learn via receiving feedback. According to \cite{yorke}, formative assessment, such as feedback, is ``vitally important to student learning. It is fundamentally a collaborative act between staff and student whose primary purpose is to enhance the capability of the latter to the fullest extent possible" (p.~496). This, and the recommendation of \cite{greenhow} to incorporate feedback in order to successfully introduce an e-learning system, supported my choice to implement automatic feedback which was available to students once a quiz deadline had passed. It would first display the question and selected answer choice, stating whether or not the selected answer choice was correct, before giving specific feedback depending on the selected answer choice. This feedback would include precise references to definitions, results or pages in the lecture notes. An example of feedback is below. This is the feedback corresponding to Example \ref{relations}. The correct answer choice for this Example is c). References to definitions etc. have been omitted for the purpose of simplicity. 

\begin{example}
\label{feedback}
$  $
\begin{enumerate}[label=\alph*)]
\item
Recall that if $ 3\mid -x+y, $ then there exists $ m\in\Z $ such that $ -x+y=3m $ by definition. Since $ m\in\Z, $ it is possible that $ m\le 0. $ Hence, our conclusion that there exists $ n\in\N $ such that $ -x+y=3n $ is incorrect. For example, consider $ x=11 $ and $ y=5. $ Then $ -x+y=-11+5=-6=2(-3)=2n $ where $ n=-3. $

Note that the proof can be easily corrected by changing $ \N $ to $ \Z ! $ This is because it is indeed correct to argue that $ -x+y\equiv 8x-5y\equiv 0 \text{ (mod } 3). $ This follows since if $ 8 \equiv -1 $ (mod $3) $ and $ x\equiv x $ (mod$ 3 $), then $ 8x\equiv -x $ (mod $3$). Similarly, if  $ 5 \equiv -1 $ (mod $3)$ and $ y\equiv y $ (mod $ 3 $), then  $ 5y\equiv -y $ (mod $3$). We can then conclude that if $ 8x\equiv -x $ (mod $3$) and  $ 5y\equiv -y, $ then $ -x+y\equiv 8x-5y \equiv 0$ (mod $ 3 $). All following conclusions are therefore correct. 
\item
In this proof, we never defined $ x $ to be an integer. Therefore, we cannot claim that if $ 8x-5x=3x $ then $ 3\mid 8x-5x. $ This is because the statement would be true if and only if $ x\in\Z. $

Note that the proof can be easily corrected by defining $ x $ to be an integer!
\item
Recall that a relation $ \sim $ is transitive if for all $ x,y,z\in\Z $ we have that if $ x\sim y $ and $ y\sim z $ then $ x\sim z. $ So we start by assuming that  $ x,y,z\in\Z $ are such that $ x\sim y $ and $ y\sim z, $ and we need to show that $ x\sim z. $

If $ x\sim y, $ then $ 8x\equiv 5y $ (mod $ 3 $) and if $ y\sim z $ then $ 8y\equiv 5z $ (mod $ 3 $). It follows that $ 3\mid 8x-5y $ and $ 3\mid 8y-5z $ which means that there exist $ m,n\in\Z $ such that $ 8x-5y=3m $ and $ 8y-5z=3n. $ We then show that $ 8x-5z=3(m+n-y). $ Since $ m+n-y\in\Z $ by assumption, we have that $ 3\mid 8x-5z. $ It then follows that $ 8x\equiv 5z $ (mod $ 3 $) which means that $ x\sim z. $
\item
The proof shows that for all $ x\in\Z $ we have that $ x\sim x. $ By definition, this shows that $ \sim $ is reflexive, so it does not show that $ \sim $ is symmetric. Therefore, our conclusion that $ \sim $ is symmetric is incorrect. 
\end{enumerate}
\end{example}

\section{Analysis and Key Results}

According to \cite{harvey} ``students are important stakeholders in the quality monitoring and assessment processes and it is important to obtain their views" (p.~19). Moreover, he says that ``feedback on specific modules or units of study provide an important element of continuous improvement" (p.~14). Therefore, in order to evaluate the effectiveness of the online quizzes, I invited students to complete a questionnaire at the end of the module in my first and second year of using the online quizzes. 

For each questionnaire, students had to select statements which they felt applied as well as answering some open-ended questions. The response rate for the first questionnaire was 59\% while the second questionnaire was completed by 37\%. Overall, I received positive and supportive comments. These included statements such as ``The quizzes give an incentive to engage with the lecture notes. This makes things more interesting.", ``The quizzes help revise stuff without putting too much pressure on me", ``Quizzes are well made" and ``I really enjoy them."

\captionsetup{width=119mm}

\begin{figure}[h!]
\centering
  \includegraphics[width=119mm]{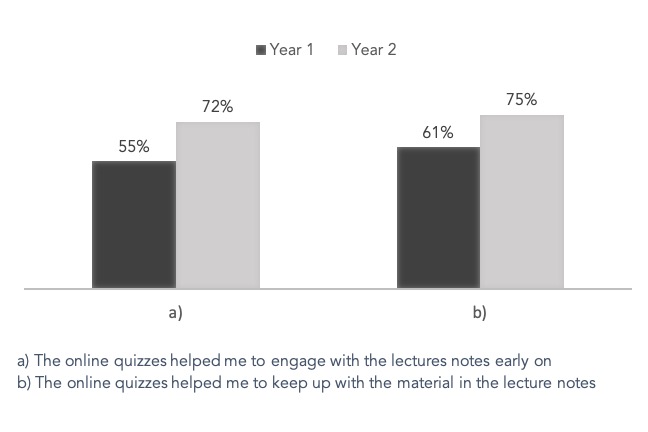}
  \caption{\, Student engagement with lecture notes}
  \label{question1}
\end{figure}

\begin{figure}[h!]
  \centering
  \includegraphics[width=119mm]{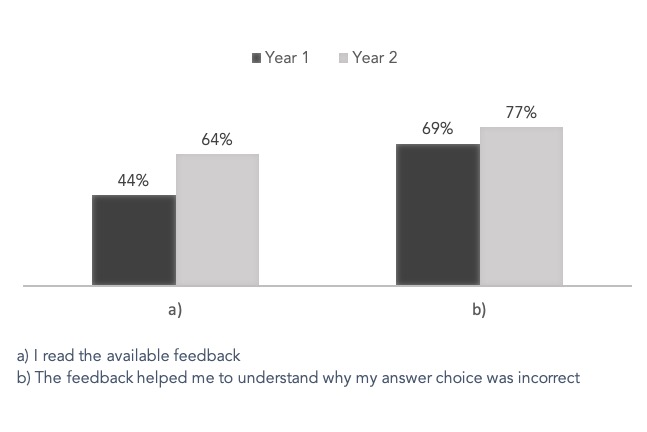}
  \caption{\, Student engagement with feedback}
  \label{question2}
\end{figure}

\begin{figure}[h!]
\centering
  \centering
  \includegraphics[width=119mm]{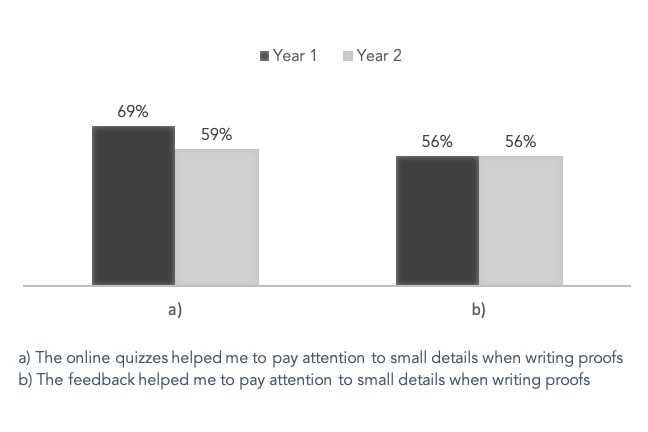}
  \caption{\, Emphasis of the small details within proofs}
  \label{question3}
\end{figure}

\begin{figure}[h!]
  \centering
  \includegraphics[width=119mm]{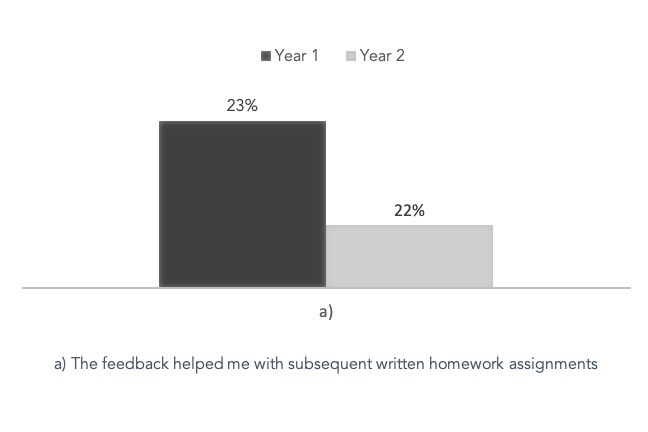}
  \caption{\, Subsequent writing of students' own proofs}
  \label{question4}
\end{figure}

First, we look at the data concerning students' engagement with their learning material, shown in Figure \ref{question1} and Figure \ref{question2}. Note that \textit{Year 1} refers to the first year of using online quizzes and \textit{Year 2} refers to the second year of using online quizzes. Figure \ref{question1} shows that 55\% of students in \textit{Year 1} found the quizzes helped them to engage with the lecture notes early on. In addition, 61\% felt that the online quizzes helped them to keep up with the material in the lecture notes. In \textit{Year 2}, the percentages increased to 72\% and 75\%, respectively. Furthermore, Figure \ref{question2} shows that while only 44\% of students in \textit{Year 1} engaged with available feedback for online quizzes, this increased to 64\% of students in \textit{Year 2}. Out of all the students who had read the feedback, 68\% of students in \textit{Year 1} thought that the feedback helped them to understand why their answer choice was incorrect. This increased to 77\% in \textit{Year 2}.

Second, we consider the data with respect to the quizzes' emphasis of the small details within proofs. Figure \ref{question2} shows that 69\% of students in \textit{Year 1} thought that the online quizzes helped them to pay attention to the small details when writing proofs. This decreased to 59\% in \textit{Year 2}. Moreover, 56\% of students, who had read the available feedback, claimed that the feedback helped them to pay attention to the small details when writing proofs. Some additional student comments on this topic included statements such as ``if you asked me to solve [the question] myself it would be fine, its not really testing my maths skill, rather my ability to see small changes in the answers". Another comment which focused on small details and written homework was that ``You can do well in the quiz without necessarily understanding the content that well. It's quite easy to spot mistakes but much harder to write a whole correct proof from scratch."

Third, we analyse the data with respect to subsequent written assignments. Here, we have that 23\% of students in \textit{Year 1} and 22\% of students in \textit{Year 2} stated that the feedback helped them with subsequent written homework assignments. An additional student comment worth mentioning was that ``I also have to write my own proof in written assignments, which I feel is more useful, though having to identify correct proofs in the online quizzes has probably helped me with that." 

Next, we have a closer look at students' non-engagement with feedback. While 27\% of students in \textit{Year 1} claimed that they were not aware that there was feedback available, this percentage decreased to 15\% in \textit{Year 2}. However, this still left 29\% of students in \textit{Year 1} choosing to not read feedback. This prompted me in \textit{Year 2} to ask students, who were aware of feedback but had chosen to not read it, for their reasons via an open-ended question. Students' responses clearly showed three trends: not needing to read feedback, time issues or timeliness of feedback. Many students claimed they had done well enough on the quiz to not need to read the feedback. Here, comments included statements such as ``I achieved high enough and understood the content in the quizzes [so] that the feedback was no longer relevant" and ``I was happy with my mark so didn't think to look at it". Comments concerning time issues included statements such as ``I did not have time so was planning to do it as revision in the holidays" and ``Haven't gotten round to it yet" and ``Intend to look back later". Last but not least, comments regarding timeliness of feedback included responses such as ``By the time it's arrived I've moved on" and ``Not available when I want it, and by the time it comes out I have more work to do so I forget''.

At the end of each questionnaire, students could leave additional comments via an open-ended question. Students' responses showed different trends in each year. Comments in \textit{Year 1} focused either on the deadline or on guessing and cheating. Comments concerning the deadline included statements such as ``It may just be nice to have a bit more time to complete the quizzes, seen as they count and with probability/analysis assessed homework, we always have a whole week." and ``Extended time to complete quizzes." Comments with respect to guessing and cheating included statements such as ``Because the questions in the online quizzes were multiple-choice and the answers were often very similar I found that it was more a case of just being able to pick one rather than actually knowing the content. "and ``Easy for people to cheat". 

Before we consider additional comments for \textit{Year 2}, I would like to briefly mention that in \textit{Year 1} I was the only lecturer using online quizzes for undergraduate  mathematics at my institution. The overwhelmingly positive feedback and the fact that 71\% of my students recommended the use of online quizzes in other modules prompted the introduction of assessed online quizzes in all but one first-year mathematics undergraduate module at my institution in \textit{Year 2}. So in \textit{Year 2}, student comments often compared the online quizzes in my module to those in other modules. They either focused on difficulty or on partial credit. Comments concerning difficulty often compared the difficulty of the quizzes with that of other modules. Such comments included statements such as ``The proofs quizzes are a lot harder than quizzes in other modules and there are a lot of questions that take a long time to do" and ``They are too hard. All the other modules are relatively easier". Comments concerning partial credit included ``They are very difficult. If you miss one thing you've lost lots of marks", ``Not a fan of online quizzes for proofs as it is very harsh with marking. If the proof selected is almost correct with minor errors all marks are still lost and this is quite frustrating" and ``I don’t like how we have big long questions and if you get one small thing wrong you lose all the marks. This does not happen in algebra."

\section{Discussion}

\subsection{Student engagement with learning material}

The data presented in the previous chapter shows how the online quizzes affected student engagement with learning material. We pay particular attention to the increased percentage of student engagement in \textit{Year 2}. Recall that some students in \textit{Year 1} claimed that answers on the quizzes were easy to guess without actually understanding the content. This encouraged me to slightly modify the online quizzes in \textit{Year 2}, changing some of the multiple-choice questions to fill-in-the-blank and matching questions. While I only changed 10\% of the questions, interestingly enough, not a single student commented on being able to guess answers in \textit{Year 2}. This implies that there was a reduced effect of guessing which further implies that, rather than being able to guess or reverse check answers, students would have to engage with the content of the questions. In order for students to engage with the content, they would have to engage with their learning material. Therefore, I would like to conclude that a reduced effect of guessing on the online quizzes in \textit{Year 2} may have contributed to the increased percentage of students engaging with their learning material. 

Moreover, I believe that there may be a correlation between student engagement with learning material and students appearing less overwhelmed by the transition from A-level mathematics to university-level mathematics. When comparing to the year prior to using online quizzes, students rarely approached me, or their tutors, with worries regarding their transition. Overall, they appeared to be less stressed with respect to the transition. Recall that according to \cite{anthony}, the study for tests and exams is a perceived key factor, by lecturers and students, for the successful transition to university-level mathematics. Therefore, I would like to conclude that there may be a correlation between students appearing less overwhelmed by the transition and the introduction of online quizzes.

Another idea worth exploring in the future which may positively affect student engagement with learning material would be to introduce one particular proof for each of the harder questions, rather than presenting four proofs, and asking several questions about that particular proof. Such proof comprehension questions could be of multiple-choice or matching question format and follow the idea of the proof understanding baseline check-list proposed by \cite{bickerton}. Questions could include, for example, asking about the type of proof being used, asking about the hypothesis of the proof, asking what definition or result is being used or asking if a justification for a particular statement is correct and providing a justification for why it is or is not. Answer choices could include, for example, statements such as ``The statement is correct because ... " or ``The statement is not correct because ... ".

\subsection{Emphasis of the small details within proofs}

Not only did the online quizzes aim at engaging students with their learning material in order to support them in their transition from A-level mathematics to university-level mathematics, but they also aimed at emphasizing the small details within proofs. The data presented in the previous chapter confirms that both of these goals were successfully achieved. Recall that 69\% of students in \textit{Year 1} and 59\% of students in \textit{Year 2} thought that the online quizzes helped to pay attention to the small details when writing proofs. Note that the statement focused on the small details when \textit{writing} proofs rather than how the online quizzes helped to pay attention to the small details within proofs in general -- a subtle difference but nevertheless a difference. In hindsight, I believe that it would have been useful to have two statements on the questionnaire concerning the small details within proofs: one which would have stated that quizzes emphasized the small details within proofs and the second one which would have stated that quizzes helped to pay attention to the small details when \textit{writing} proofs. This would have been useful to establish a possible gap between students paying attention to the small details when \textit{reading} proofs and being able to apply this knowledge when \textit{writing} proofs. Recall that proof validation has been shown to have a positive impact on students' ability to construct their own proofs \cite[]{powers}, so the possibility of a gap is supported by the data which shows that only a low percentage of students thought that the online quizzes helped them with subsequent written homework assignments or online quizzes. It is further evidenced by the fact that while 69\% of students in \textit{Year 1} and 77\% of students in \textit{Year 2} indicated that the feedback helped them to understand why their particular answer choice was incorrect, only 56\% of students in both years indicated that the feedback helped them to pay attention to the small details when writing proofs. Moreover, student comments such as ``It's quite easy to spot mistakes but much harder to write a whole correct proof from scratch" further support the possibility of a gap. This leads me to conclude that while the quizzes successfully achieved their goal of emphasizing the small details within proofs, I believe it to be beneficial to try and close the gap so that students may learn how to apply the knowledge they learnt when writing their own proofs. While I already used examples in problem classes similar to questions which would appear in the upcoming quiz, I believe it would be beneficial to also apply this principle retrospectively: to point out the link between previous quiz questions and, for example, exercises in problem classes. 

I believe that the automatic feedback given to students further emphasized the small details within proofs by explaining why a particular answer choice was correct or incorrect. While only 59\% of students, who had read the available feedback, claimed that the quizzes helped them to understand the small details when \textit{writing} proofs, it is clear from their statements that they realised that small details can have a big impact, that is, they can be a cause of losing marks. Since an answer choice of a multiple-choice question is either correct or incorrect, this question type especially emphasizes how an entire proof may collapse due to one small detail. While students did not seem particularly satisfied with this impact, I believe that it nevertheless successfully emphasized the importance of the small details within proofs. And interestingly, when I changed some of the multiple-choice questions to either fill-in-the-blank or arrange-a-proof-in-order matching questions in order to reduce the effect of guessing in \textit{Year 2}, the percentage of students claiming that the online quizzes helped to emphasize the small details within proofs decreased, so there may be a possible correlation. Because of this, I think it to be useful to not give partial credit on multiple-choice questions. However, upon reflection, I think it to be perfectly reasonable to allow partial credit for multiple-answer and matching questions. 

\subsection{Student engagement with feedback}

According to \cite{gedye}, students rarely engage with feedback despite its benefits, so perhaps the fact that only 44\% of students in \textit{Year 1} read the available feedback was to be expected. However, the fact that 27\% of students claimed that they were not even aware that there was feedback available is still surprising for several reasons. First of all, I announced the availability of feedback on several occasions in problem classes in \textit{Year 1}. I showed students how to find this feedback on Blackboard, and these instructions were recorded via Mediasite so they were available via Blackboard. Second of all, I added a PDF file to Blackboard with instructions on where and when to find feedback. This PDF file was available on the same page as the online quizzes. The fact that 56\% of students in \textit{Year 1} did not engage in feedback prompted me in \textit{Year 2} to repeatedly mentioned the availability of feedback in problem classes as well as lectures. In addition, I went through the instructions on how to view feedback on Blackboard on at least a couple of occasions. I believe that this repeated emphasis likely contributed to the 20\% increase in the percentage of students reading feedback and the 12\% decrease in the percentage of students not being aware that there was feedback available.

Another idea worth exploring in the future which may increase student engagement with feedback would be to release a copy of the online quizzes at the end of the term for revision. Students would then be able to re-attempt the quizzes (without them counting towards a mark) and receive immediate feedback. While this may not necessarily increase student engagement with feedback during term time, it may well increase overall engagement. 

\subsection{Timeliness of feedback}

Recall that some students claimed that they had not read the available feedback because it was not accessible when it was relevant. Further, recall that the feedback was automatically available to students once a quiz deadline had passed. In \textit{Year 1}, this deadline was 48 hours, but considering the fact that several students commented that they would like to have a longer deadline and considering the fact that students, in general, would have a one week deadline for written coursework, I only thought it to be fair to introduce a one week deadline in \textit{Year 2}. While I no longer received any concerns regarding the deadline of the quizzes, this extension unfortunately meant that if students submitted early within the one week deadline they had to wait until the deadline had passed to be able to read the feedback. \cite{irons} suggests that timeliness of feedback is a key feature in students' learning process, but according to \cite{bayerlein} ``students do not distinguish between timely feedback and extremely timely feedback" (p.~916). Further, Bayerlein says that ``undergraduate students do not recognise any improvements in timeliness when feedback is provided in less than 11.5 days" (p.~924). This critical time frame of two weeks was confirmed by \cite{bohnacker} and \cite{brown}. When we additionally consider the fact that students receive feedback on formative written homework assignments after one week and on summative written coursework assignments within three weeks, whereas the maximum amount of time to receive feedback on the online tests was one week, then I feel the current timeliness of quiz feedback to be fair and constructive.

\subsection{Level of difficulty}

Last but not least, I would like to address student comments concerning the difficulty of the quizzes. While I did not receive any such comments in \textit{Year 1}, a few students in \textit{Year 2} claimed that the quizzes were too difficult. Since I only changed 10\% of the quiz questions, it does not seem plausible that this change could account for such a perceived increase in the level of difficulty. However, considering the fact that in \textit{Year 1} I was the only module leader to use online quizzes and that this prompted all but one first-year undergraduate mathematics module to introduce online quizzes in \textit{Year 2} combined with the fact that students commented on quizzes being more difficult than in other modules, I think it is safe to conclude that students found the online quizzes in my module to be more difficult when being compared to online quizzes in other modules. I would also like to add that while students commented on the perceived difficulty of the quizzes, the overall average of student marks was not affected compared to that of students in \textit{Year 1}. When further comparing average students marks from the year(s) when using online quizzes to the year(s) prior to using online quizzes, then, again, we do not observe a noticeable difference. While this does not necessarily disagree with studies which showed an increase in summative exam scores when using e-assessments \cite[]{angus, dobson, kibble, salas}, it also does not necessarily support the argument either. However, it is debatable whether the comparison of average student marks is, in this case, particularly conclusive since, for example, the cohort, the content of the exam, some of the quiz questions, etc. varied each year. Moreover, the aforementioned studies focused on subject areas other than mathematics, and as previously established, mathematics is a discipline with unique features when it comes to assessment and teaching \cite[]{lms}.

\section{Conclusion}

In conclusion, I believe that the online quizzes successfully achieved their intended goals of, firstly, providing students with an extrinsic motivation to engage with learning material early on in order to academically support them with their transition from A-Level mathematics to university-level pure mathematics and, secondly, of emphasizing the importance of the small details within proofs. While only up to 23\% of students seemed to think that the online quizzes helped them with the subsequent writing of their own proofs, I believe that this is due to a gap between learnt knowledge and knowing how to apply this learnt knowledge. However, I do believe that, in the future, it will be possible to easily bridge that gap by, for example, linking learnt knowledge from the online quizzes with exercises in problem classes. 

Furthermore, while I do not believe that Blackboard is flawlessly adapted to mathematics and that it would certainly benefit the integration of STACK, I think that its use for online quizzes is sustainable until more well-adapted e-learning systems for mathematics are developed. Advantages of using Blackboard are that it supports the basic use of \LaTeX, automatic marking which is especially great for large cohorts, a variety of questions types, automatic feedback tailored to students' selected answer choices, and automatic randomisation of questions and answer choices so that quizzes may be re-used in subsequent academic years.

Henceforth, I believe that institutions, which only have one summative end-of-term exam, would benefit from using online quizzes in first-year undergraduate pure mathematics courses -- supplementary to formative written homework assignments -- as a means to engage students with learning material early on and to support them in their proof learning process. For the successful introduction of online quizzes, I would like to add several recommendation. The final four recommendations are based on my own findings and personal observations:

\begin{itemize}
\item[$ \circ $]
to encode common students' misconceptions \cite[]{greenhow}
\item[$ \circ $]
to schedule quizzes periodically to allow for continuous engagement \cite[]{savander}
\item[$ \circ $]
to provide detailed feedback tailored to selected answer choices \cite[]{greenhow}
\item[$ \circ $]
to randomise questions and answer choices \cite[]{greenhow}
\item[$ \circ $]
to make quizzes summative by counting towards a small percentage of students' mark to provide an extrinsic motivation for students to engage
\item[$ \circ $]
to use open-book format
\item[$ \circ $]
to link learnt knowledge with, for example, written homework assignments or exercises in problem classes or lectures
\item[$ \circ $]
to allow for partial credit whenever possible
\end{itemize}

\section*{Acknowledgements}

This project was supported by the Bristol Institute for Learning and Teaching (BILT) Associateship scheme. I would also like to extend my gratitude to Dr. Paola Iannone, Senior Lecturer at Loughborough University, for sharing her pearls of wisdom on an earlier version of the manuscript.

\bibliography{ref}
\bigskip

\appendix

\section*{APPENDIX}
\label{questionsamples}

Below are some more examples of the different question types used for the online quizzes. Example \ref{multiplechoice} is an example of a multiple-choice question and Example \ref{multipleanswer} is an example of a multiple-answer question.

\begin{example}
\label{multiplechoice}
Let $X=\{a,b\}$ and $Y=\{1,2\}$. What is the Cartesian Product of $X$ and $Y$?
\begin{enumerate}[label=\alph*)]
\item
$X\times Y=\{(a,1),(b,2)\}$
\item
$X\times Y=(a\cdot 1, b\cdot 2)$
\item
$X\times Y=\{(a,1),(a,2),(b,1),(b,2)\}$
\item
$X\times Y= a\cdot 1 +b\cdot 2$
\end{enumerate}
\end{example}

\begin{example}
\label{multipleanswer}
Denote by $ \mathcal{P}(A) $ the power set of $ A. $ Which of the following sets are countable? Select \textbf{all} that apply. 
\begin{enumerate}[label=\alph*)]
\item 
$\mathcal{P}(\mathbb{Q})$
\item
$\mathbb{R}_+$
\item
$\left\{ x\in\mathbb{R}: x=\frac{n^2}{\pi}, n\in\N \right\}$
\item
$ \N^n $ for $ n\in\N $
\item
$ \mathcal{P}(\mathcal{P}(\{2,3\})) $
\end{enumerate}
\end{example}

\end{document}